\documentclass[english]{article}
\usepackage[T1]{fontenc}
\usepackage[latin9]{inputenc}
\usepackage{geometry}
\usepackage{float}
\usepackage{amsmath}
\usepackage{amssymb}
\usepackage{graphicx}

\makeatletter

\providecommand{\tabularnewline}{\\}

\numberwithin{equation}{section}
\newcommand{\lyxaddress}[1]{
\par {\raggedright #1
\vspace{1.4em}
\noindent\par}
}

\@ifundefined{date}{}{\date{}}
\usepackage{url}

\makeatother

\usepackage{babel}
\begin{document}

\title{A Two-Level Variant of Additive Schwarz Preconditioning for Use in
Reservoir Simulation}

\author{Haran Jackson%
\thanks{Current address: Hertford College, Catte Street, Oxford, OX1 3BW%
}, Michele Taroni and David K. Ponting%
\thanks{Corresponding author: david.ponting@emerson.com%
}}

\maketitle

\lyxaddress{\begin{center}
Roxar Ltd, Emerson Process Management, Northbrook House, Oxford Science
Park, Oxford, OX4 4GA
\par\end{center}}
\begin{abstract}
The computation time for reservoir simulation is dominated by the
linear solver. The sets of linear equations which arise in reservoir
simulation have two distinctive features: the problems are usually
highly anisotropic, with a dominant vertical flow direction, and the
commonly used fully implicit method requires a simultaneous solution
for pressure and saturation or molar concentration variables. These
variables behave quite differently, with the pressure feeling long-range
effects while the saturations vary locally. In this paper we review
preconditioned iterative methods used for solving the linear system
equations in reservoir simulation and their parallelisation. We then
propose a variant of the classical additive Schwarz preconditioner
designed to achieve better results on a large number of processors
and discuss some directions for future research.
\end{abstract}

\section{Introduction}

In terms of computing time, the primary element of the majority of
reservoir simulations is the solution of a large sparse set of non-symmetric
linear equations. The most common mode for reservoir simulation is
three-phase fully implicit. Over a given time step in the simulation,
a set of non-linear conservation equations is solved, the solution
variables typically being cell pressures and either saturations or
molar densities (see, for example, \cite{Chen07}). It is also possible
to treat just the pressures implicitly, replacing some non-linear
functions such as fluid densities and capillary pressures by their
start of time step values. Such a method, usually referred to as IMPES
(implicit pressure explicit saturation), is now mainly used for compositional
reservoir simulation, in which multiple components are tracked using
an equation of state to characterise the fluid. IMPES methods are
given to oscillatory instability, and fully implicit methods are generally
preferred by users. 

The required system of conservation equations is typically derived
from a finite volume discretisation, and consists of mass accumulation,
flow, and well injection and production terms. The system is conservative,
and is commonly expressed in residual form, $R(X)=0$, where $X$
is the set of solution variables to be solved for. In solving these
equations using a Newton-Raphson method, a set of solution changes
$x$ is obtained by solving the set of linear equations given by

\begin{equation}
x=-J^{-1}R\left(X\right)
\end{equation}

where

\begin{equation}
J=\nabla R\left(X\right)
\end{equation}

$J$ is the Jacobian. It is common to treat the set of solution variables
associated with a given reservoir simulation cell (for example the
pressure, oil saturation and gas saturation) as a single strongly
coupled sub-vector. So for a three-phase black oil problem the elements
of $X$ and $R$ would be 3-component sub-vectors and the elements
of $J$ would be 9-component sub-matrices. Given the large size of
the ensuing linear system (hundreds of thousands of cells are common
in industry, millions becoming so) direct methods are not feasible,
and so iterative methods are used.

In this paper we begin in Section 2 by reviewing the iterative methods
typically used in reservoir simulation packages, before discussing
their adaptation and performance in a parallel computing environment
in Section 3. In Section 4 we then present some recent work aimed
at improving the parallel performance of multi-phase linear solution
by extending the overlap zones to cover the full problem in a coarsened
way, similar to a coarse-grid correction but ensuring that local saturation
variables are not smeared. Results are presented in Section 5, before
discussing our findings and suggesting future avenues of research
in Section 6.

\section{Review of iterative methods for solving linear system in reservoir
simulation}

The linear sub-problem can also be expressed in residual form as $r=b-Ax=0$,
with $b=R$, the current non-linear residual, and $A=-J$. In the
early days of reservoir simulation stationary iterative methods such
as SOR were common \cite{Wat71}, but soon preconditioned conjugate
gradients became standard, typically accelerated using ORTHOMIN \cite{Vin76},
GMRES \cite{Saad86}, or BiCGStab \cite{vdV92}. Preconditioning is
usually performed by selecting an invertible approximation to $A$,
$B\sim A$, such that the inversion of $B$ is reasonably practical.
The most common choice of $B$ is an incomplete LU factorisation with
some degree of fill-in, often zero \cite{Saad03}.

In the context of fully implicit reservoir simulation it is typical
for ILU methods to be applied to the block structure of the matrix,
with the dense sub-matrices inverted directly. Furthermore, the issue
of material balance is important in reservoir simulation. In light
of this, it is possible to choose $B$ such that the column sum of
the error matrix $E=B-A$ is zero; this improves convergence and corresponds
to selecting a solution in which $\sum r=0$. This has a useful physical
interpretation in that it corresponds to zero mass conservation error.
As a reservoir simulation can involve many thousand steps, it is useful
if the mass accumulation error superconverges: in some reservoir simulation
formulations the mass balance error can be zero after one non-linear
iteration if using a linear solver which zeros the linear residual
sum.

\begin{figure}
\centering{}\includegraphics[scale=0.5]{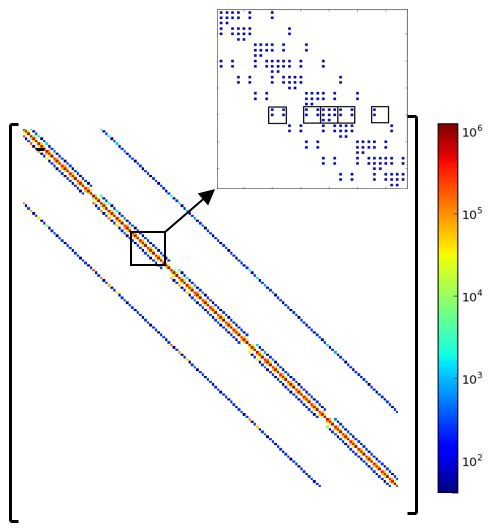}\caption{Sparsity pattern of a Jacobian from a typical fully implicit problem.
Colour scale indicates size of $3\times3$ implicit blocks.}
\end{figure}

The Jacobian matrices arising from reservoir simulation problems are
diagonally dominant, but marginally so, with strong bands representing
flow between neighbouring cells \textendash{} an example is shown
in Figure 1. Furthermore, an oil reservoir may be tens of kilometres
in areal extent, but a few hundred metres in vertical extent. Flow
in the vertical direction is over a shorter distance than in the areal,
and occurs over a larger area. The result is that off-diagonal matrix
elements representing vertical flow are particularly strong. In light
of this, a popular variant of ILU within the context of reservoir
simulation is nested factorisation, a series of nested LU factorizations,
which is particularly suited to the strong asymmetry of the problem
\cite{App83}. Nested factorisation produces a preconditioning matrix
with no error terms in the vertical direction, as shown in Appendix
A, and by eliminating errors in the dominant flow direction results
in a very effective and robust preconditioning. Furthermore, by requiring
no fill-in, the algorithm is memory efficient. However, it is less
suitable for reservoirs with horizontal wells, which destroy the banded
sparsity pattern, and those with large heterogeneity, especially in
the horizontal direction.

An alternative to preconditioned iterative methods are algebraic multigrid
(AMG) methods. These construct a solution based on the application
of corrections obtained on a series of successively coarser grids.
The corrections are frequently regarded as smoothings, and constructed
using Gauss-Seidel or ILU techniques. Such methods are generally effective
on the long-range terms which arise in the solution of elliptic equation
systems. In reservoir simulation terms this corresponds to the pressure
variable part of the solution. However, multigrid methods are less
suited to the treatment of saturation or composition variables: these
exhibit sharp variations in fluid displacement fronts, and the use
of a coarse grid method introduces solution smearing. The result is
that, for fully implicit simulation, multigrid methods are usually
applied in combination with a second saturation solver: first an approximate
pressure-only-equation is set up and solved (using, for example, the
constrained pressure residual formulation \cite{Wal83}); then a solver
such as ILU(0) is used to obtain the full pressure and saturation
solution to the full system \cite{Stu07,Ham08}. It is then common
to use this two-stage solution as a preconditioner within an iterative
method.

\section{Parallelisation}

With the advent of parallel computing architectures in the last twenty
or so years reservoir simulation solvers have been adapted to work
in parallel. Unfortunately, the use of an invertible LU approximation
to the matrix, such as nested factorisation, tends to involve recursive
back-substitution techniques which are difficult to run in parallel.
One common method for improving the performance of LU factorisation
methods when running in parallel is the use of multi-colour domains.
Examples are simple block red-black ordering such as parallel nested
factorisation \cite{Pon96} and more complex multi-colour methods
such as the JALS solver \cite{App11}. These preserve the virtue of
nested factorisation that error terms in one direction are eliminated,
but the preconditioning becomes less effective in the other directions
as the number of processors is increased, increasing the iteration
count. In addition, the presence of features such as horizontal wells
can yield matrix elements which do not fit naturally into the multicolour
ordering. Multigrid solvers are more effective as parallel solvers,
the number of linear iterations remaining more constant as the number
of processors is increased. However, a significant amount of inter-process
communication is required to form and solve for the corrections on
the coarser grids using all the processors; and an effective parallel
method to find the full pressure and saturation solution after the
multigrid pressure correction has been applied is still required.
Overall, the effect is that solver speed-up still eventually degrades
with increasing number of processors.

Whilst the use of up to around 60 processors has been common for some
years, computers with hundreds of thousands of processors are now
becoming available, so there has been a drive towards finding a new
generation of solvers that can scale to work on large clusters of
processors. A class of methods that are attractive in this context
use sparse approximate inverse (SAI) solvers, which have been widely
investigated in the academic literature \cite{Ben00,Chow01}. These
are radically different in that they work in terms of creating approximations
to $A^{-1}$ directly rather than an invertible approximation to $A$.
In a SAI method the construction of the approximate inverse is thus
a completely parallel task; and the construction of a new search direction
is done using a matrix multiplication operation. However, in order
to render the method practical the approximate inverse must be constrained
to be sparse, which means that when constructing a new linear solver
search direction a given solution variable is only aware of a small
subset of elements of the residual $r$. Generally SAI methods are
similar in performance to ILU methods and have not been heavily used
in practical reservoir simulation. 

Current methods for solving reservoir simulation linear equations
on highly parallel systems include the LSPS solver \cite{Fung07},
which uses a power expansion of the inverse of the matrix $A$ to
obtain a better approximation to the inverse matrix $A^{-1}$:

\begin{equation}
B^{-1}=\left(I+\sum_{i=1}^{n}\left(-T^{-1}E\right)^{i}\right)T^{-1}
\end{equation}

where

\begin{equation}
A=T+E
\end{equation}

and the use of additive Schwarz methods in which each processor uses
its own block of the entire matrix, plus some elements of the matrix
from other processors by overlapping the computational domains:

\begin{equation}
B_{AS}^{-1}=\sum_{i=1}^{n}R_{i}^{T}A_{i}^{-1}R_{i}
\end{equation}

where

\begin{equation}
A_{i}=R_{i}AR_{i}^{T}
\end{equation}

and $R_{i}$ is the rectangular restriction matrix to subdomain$i$.

The draw-back with additive Schwarz methods is that they do not remove
small eigenvalues of the coefficient matrix, corresponding to low
frequency modes arising from the long range nature of pressure effects
in reservoir simulation. Two-level additive Schwarz methods (see,
for example, \cite{Smith96}) attempt to alleviate this problem by
adding a coarse-grid correction, which acts like the high level coarse
grid corrections in a multigrid method:

\begin{equation}
B_{AS2-level}^{-1}=\sum_{i=1}^{n}R_{i}^{T}A_{i}^{-1}R_{i}+R_{C}^{T}A_{C}^{-1}R_{C}
\end{equation}

where

\begin{equation}
A_{C}=R_{C}AR_{C}^{T}
\end{equation}

and $R_{C}$ is a suitable global coarse-grid operator. Alternatively,
one can apply the coarse grid correction successively, leading to
an overall \textquotedblleft{}two-stage\textquotedblright{} preconditioner

\begin{equation}
B_{AS2-stage}^{-1}=B_{AS}^{-1}+B_{C}^{-1}\left(I-AB_{AS}^{-1}\right)
\end{equation}

where

\begin{equation}
B_{C}^{-1}=R_{C}^{T}A_{C}^{-1}R_{C}
\end{equation}

\pagebreak{}

\section{Boundary conditioning}

In the context of multi-phase flow, such a coarse-grid correction
suffers from the same problem as multigrid: the saturations are smoothed
along with the pressure. In order to circumvent this problem, we propose
a modified additive Schwarz preconditioner in which the operator $R_{i}$,
rather than discarding parts of the domain outside of the ith subdomain,
coarsens them. In matrix terms the operator may be written as

\begin{equation}
R_{i}=\left(\begin{array}{ccccccc}
R_{C,1}\\
 & R_{C,2} &  &  &  & 0\\
 &  & \ddots\\
 &  &  & I_{i}\\
 &  &  &  & \ddots\\
 & 0 &  &  &  & R_{C,n-1}\\
 &  &  &  &  &  & R_{C,n}
\end{array}\right)
\end{equation}

where $R_{C,i}$ are local coarse-grid operators. In this way each
processor block sees the whole problem to some extent, but avoids
applying a coarse grid correction to the saturation variables on that
processor. An example of this process for a simple one-dimensional
decomposition of a two-dimensional domain is shown in Figure 2.

\begin{figure}[H]
\centering{}\includegraphics[scale=0.5]{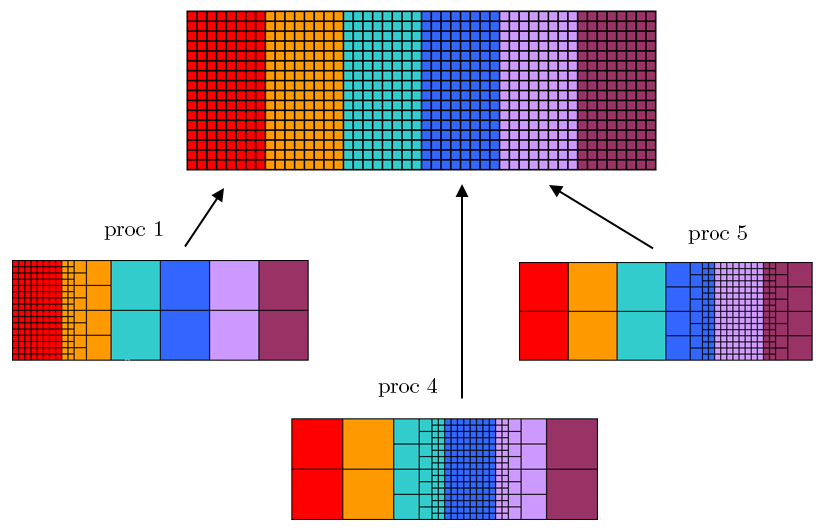}\caption{Sketch showing a six-way one-dimensional decomposition of a two-dimensional
uniform Cartesian grid (top) and how each processor sees only a coarsened
version of the entire problem in which cells outside its own subdomain
are lumped together (for clarity, only processors 1, 4, and 5 are
shown)}
\end{figure}

\newpage{}

The method proposed differs from traditional methods such as parallel
nested factorisation in not using a single invertible approximating
matrix $B$; instead each processor constructs and inverts its own
approximation to the full matrix. In this way the method is similar
to a block sparse approximate inverse; however the reduction in the
dimension of each processor matrix is done using coarsening rather
than sparsification, so that the reduced system is a coarse rather
than a sparse approximate inverse. A similar approach has been taken
by \cite{Bank02} in the context of finite element methods; however
our approach differs in that the coarsening is performed by simple
matrix element summation \textendash{} the motivation is that solving
a system formed in this way corresponds to finding a solution in which
residual sums over the coarsened regions are forced to zero, so that
material is conserved in a coarse grid sense (see Appendix B).

Once each processor has inverted its own approximation to the full
matrix the full linear search direction is obtained by projecting
out the on-processor elements for each domain, as in the restricted
additive Schwarz method \cite{Cai99}. The solution obtained for the
off-processor coarsened variables is discarded; it exists only to
set up the boundary conditions for the on-processor variables.

The method can thus be represented as:

\begin{equation}
B^{-1}=\sum_{i=1}^{n}P_{i}B_{i}^{-1}
\end{equation}

where

\begin{equation}
B_{i}^{-1}=R_{i}^{T}A_{i}^{-1}R_{i}
\end{equation}

$P_{i}$ is a projection operator into the range of variables corresponding
to a processor $i$. The full approximate inverse matrix is thus constructed
by taking the rows corresponding to each processor from each of the
different preconditionings used on each processor. In practice, the
are never explicitly formed: instead the required rows from each are
combined to form the full preconditioner. This row-wise combination
means that if the original matrix $A$ were symmetric, the preconditioning
matrix would not be. A similar effect occurs in SAI preconditioning,
in which the approximate inverse is constructed on an independent
column by column basis.

In matrix terms, at each iteration of the outer iterative scheme each
processor $i$ solves

\begin{equation}
\left(\begin{array}{cc}
A_{i,i} & A_{i,C}\\
A_{C,i} & A_{C,C}
\end{array}\right)\left(\begin{array}{c}
\Delta x_{i}\\
\Delta x_{C}
\end{array}\right)=-\left(\begin{array}{c}
\Delta r_{i}\\
\Delta r_{C}
\end{array}\right)
\end{equation}

where the subscripts $C$ indicates coarsened terms coming from other
processors; the global search direction is then taken to be

\begin{equation}
\Delta x=\sum_{i}\Delta x_{i}
\end{equation}

Like multigrid, inter-processor communication is required to form
the coarsened matrices, but this occurs only at the start of each
non-linear iteration. At each linear iteration only the updated coarsened
residuals must be communicated. Furthermore, we note that although
ideally refinement might steadily increase away from a given processor
block, to obtain reasonable computation efficiency matrix elements
must be lumped before communication. If all processors required a
different coarsening then this algorithm would not scale, and so it
is necessary to allow only a fixed number of coarsenings. In particular
we use a \textquoteleft{}far-field\textquoteright{} coarse approximation,
shared by all processors, and in addition a more refined \textquoteleft{}near-field\textquoteright{}
approximation is constructed and shared between neighbouring processors;
a sketch of the operations performed is given in Figure 3. This two-level
approach is motivated by the multi-phase physics of the problem, with
the near-field designed to capture local saturation variables while
the far-field characterises the global pressure variations. The idea
is to include long range pressure effects in an estimation of the
values of the near-field cells. These then provide reasonable boundary
values for the calculation on a particular processor.

\begin{figure}[H]
\centering{}\includegraphics[scale=0.5]{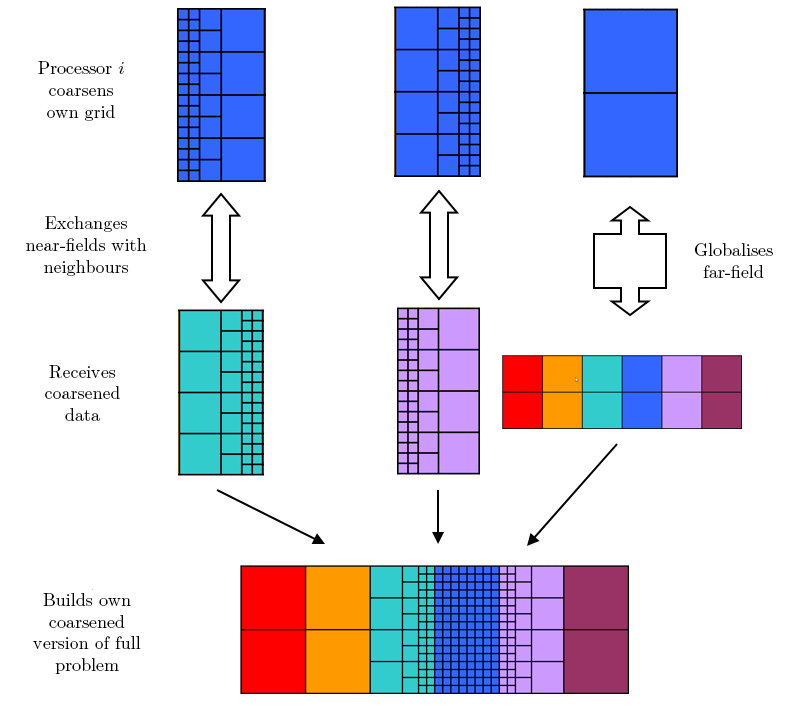}\caption{Sketch showing how each processor coarsens its own grid in three different
ways, exchanges data with other processors, and assembles its coarsened
version of the full problem}
\end{figure}

It is not a requirement that the far-field blocks correspond to processor
domains. Generally, the smaller the blocks the better the boundary
values that are supplied to each processors; but the price paid is
that the total number of variables solved for increases.

\section{Results}

In order to test the theoretical performance of the method, tests
were performed using a Python implementation of the algorithm, with
the subdomain problems $A_{i}\Delta x_{i}=-\Delta r_{i}$ solved exactly
by Gaussian elimination. The tests were performed on a number of matrices
arising from reservoir simulation applications, whose properties are
listed in Table 1. The first benchmark test case is the publicly available
ORSREG1 matrix from the Harwell-Boeing matrix collection \cite{H-B},
taking the unit vector as the right-hand side. The others are variants
of the benchmark test case of Aziz Odeh \cite{Odeh81} using water
injection instead of gas injection, with four different well configurations,
as shown in Figure 4. The original grid ($10\times10\times3$) has
been refined by a factor of two or three in each direction to give
larger matrices, although the implementation restricted us to matrices
of size less than 105. Both IMPES and fully implicit versions of the
problem were considered.

We compare preconditioning using just a near-field (which is effectively
an additive Schwarz method), just a far field, (which is a coarse
grid correction, but applied within a single stage preconditioning)
and a combination of the near and far field methods. For simplicity,
the domain decomposition is done only in the $x$-direction.

\begin{table}[H]
\centering{}%
\begin{tabular}{|c|c|c|c|c|}
\hline 
Matrix Name & Type & Size & Number of Non-Zeros & Condition Number\tabularnewline
\hline 
\hline 
ORSEG1 & IMPES & 2205 & 14133 & 6745\tabularnewline
\hline 
SPE1a & IMPES & 8100 & 53874 & 349422\tabularnewline
\hline 
SPE1b & IMPES & 8100 & 54430 & 52436\tabularnewline
\hline 
SPE1c & IMPES & 8100 & 54432 & 379882\tabularnewline
\hline 
SPE1d & IMPES & 8100 & 54690 & 43969\tabularnewline
\hline 
SPE1a\_F & FULLI & 9600 & 102878 & 226904602\tabularnewline
\hline 
SPE1b\_F & FULLI & 9600 & 108964 & 38324890\tabularnewline
\hline 
SPE1c\_F & FULLI & 9600 & 107790 & 44489039\tabularnewline
\hline 
SPE1d\_F & FULLI & 9600 & 129399 & 5941993\tabularnewline
\hline 
\end{tabular}\caption{Properties of test matrices used}
\end{table}

\begin{figure}[H]
\centering{}\includegraphics[scale=0.5]{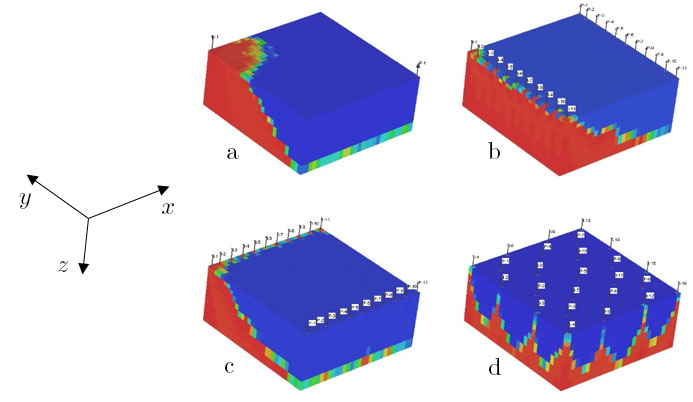}\caption{Variants of SPE1 benchmark from which matrices SPE1a,b,c,d are taken.
The colour represents the water saturation ranging from zero (blue)
to one (red). In all cases the domain decomposition was performed
in the $x$-direction.}
\end{figure}

In Figure 5 the linear iteration count for solving the ORSREG1 matrix
is shown. $K$ is the number of processors used, $C$ is the size
of the coarsening (e.g. $C=2$ means using coarsening blocks of $2\times2\times2$)
and \textquotedblleft{}total blocks\textquotedblright{} refers to
the average size of each $A_{i}$ and is a measure of how much work
goes into building the pre-conditioner. As can be seen, the combination
of near and far field consistently yields the minimum number of iterations.
However, this does involve solving for more additional variables than
the other two methods, which offsets the gain in iteration count to
some extent. Nevertheless, the solver compares favourably to the classical
two-level Schwarz method and alternative two-stage method, as shown
in Figure 6.

\begin{figure}[H]
\centering{}\includegraphics[scale=0.95]{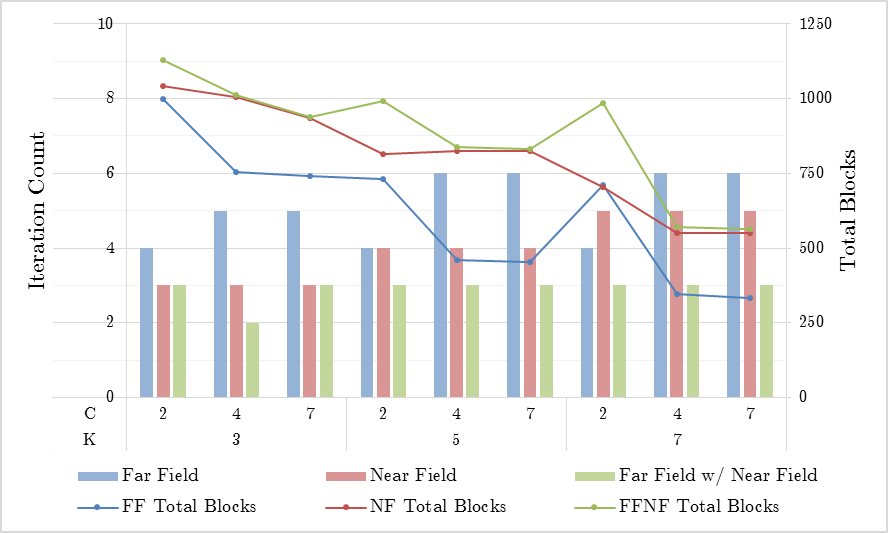}\caption{Iteration count (histogram values) against number of processors $K$
and size of coarsening $C$ for the ORSREG1 benchmark. The average
size of the subdomain problems (\textquotedblleft{}total blocks\textquotedblright{})
is also shown.}
\end{figure}

\begin{figure}[H]

\centering{}\includegraphics[scale=0.95]{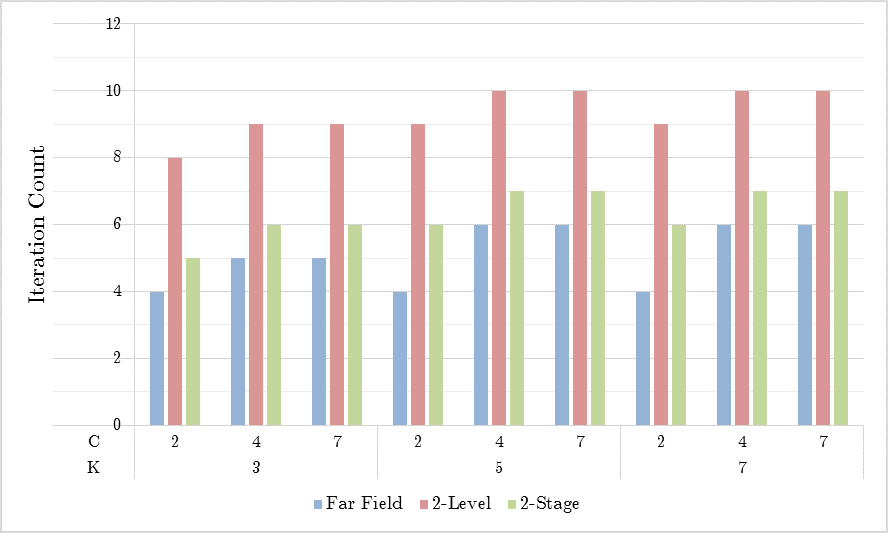}\caption{Iteration count against number of processors and size of coarsening
for the ORSREG1 benchmark, comparing the far-field method to the two-level
and two-stage additive Schwartz preconditioners}
\end{figure}

For simplicity we have shown results without overlap, so as to compare
including coarsened information in our far-field method compared to
the classical methods \textendash{} however similar results are obtained
when an overlap (or near-field) is added. The same behaviour was seen
for all of the matrices tested, with the near+far method consistently
giving the best iteration count \textendash{} examples are shown in
Figures 7 and 8.

\begin{figure}[H]
\centering{}\includegraphics[scale=0.95]{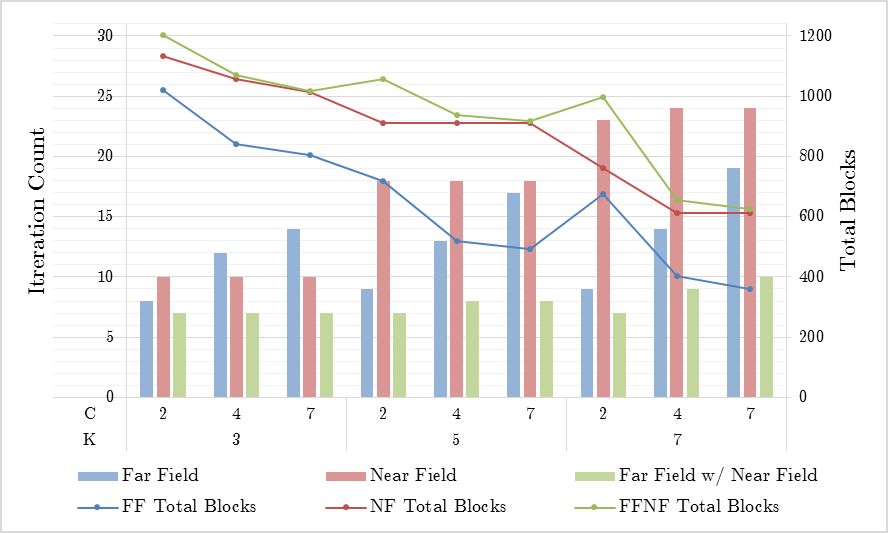}\caption{Iteration count (histogram values) against number of processors $K$
and size of coarsening $C$ for the fully implicit SPE1a\_F benchmark
refined to $20\times20\times6$ cells. The average size of the subdomain
problems (\textquotedblleft{}total blocks\textquotedblright{}) is
also shown.}
\end{figure}

\begin{figure}[H]
\centering{}\includegraphics[scale=0.95]{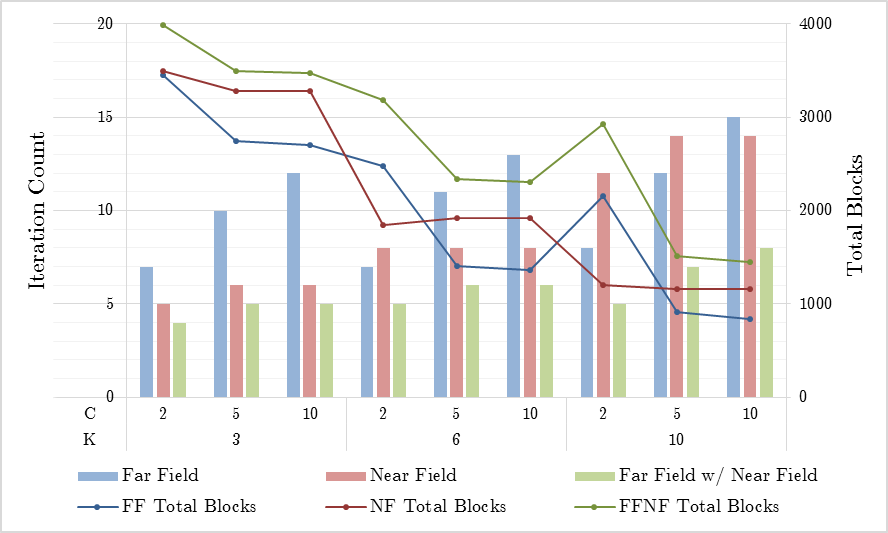}\caption{Iteration count (histogram values) against number of processors $K$
and size of coarsening $C$ for the IMPES SPE1d benchmark refined
to $30\times30\times9$ cells. The average size of the subdomain problems
(\textquotedblleft{}total blocks\textquotedblright{}) is also shown.}
\end{figure}

\newpage{}

Finally, tests on a two-dimensional decomposition (in both $x$ and
$y$ direction) showed this to be superior to those based on a one-dimensional
split, see Figure 9. The only case in which the one-dimensional split
worked as well as the two-dimensional split was case b, where the
primary direction in which the pressures and saturation vary is orthogonal
to the direction in which the problem is divided into processors.

\begin{figure}[H]

\centering{}\includegraphics[scale=0.95]{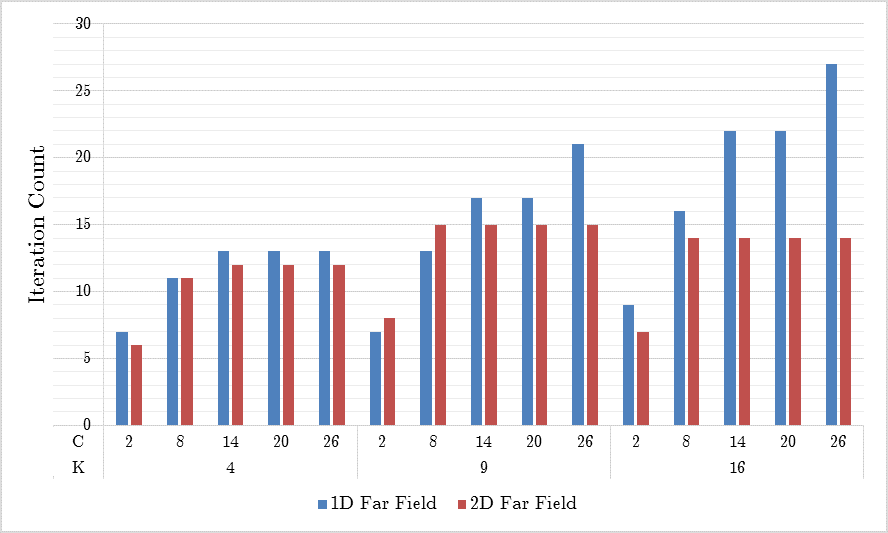}\caption{Iteration count against number of processors $K$ and size of coarsening
$C$ for the IMPES SPE1d benchmark refined to $30\times30\times9$
cells. Results are using the far-field solver with domain decomposition
in one and two dimensions.}
\end{figure}

\section{Discussion}

The method proposed is an extension of additive Schwarz which takes
account not only of the neighbouring processor\textquoteright{}s cells,
but also of cells further away in a coarsened sense. As the coarse
grid problem is effectively solved in a slightly different form on
all the processors, this background far field needs to be fairly coarse,
both to minimise the number of off-processor cells and the communication
time required to pass this information to all the processors. The
advantage is that the whole problem is seen, in some sense, by all
the processors. The improvement in preconditioning is offset by the
number of \textquoteleft{}off-processor\textquoteright{} variables
which must be treated by each processor along with the \textquoteleft{}on-processor\textquoteright{}
ones that are actually used. Clearly, if the problem is large, there
is a case for a second deeper level of coarsening of the far-field
in distant regions far from a given processor. 

The results presented in this paper were all on small problems with
less than 105 variables, with the subdomain problems $A_{i}\Delta x_{i}=-\Delta r_{i}$
solved exactly by Gaussian elimination. This gives us the optimal
iteration count for the proposed method, focusing solely on the effect
of boundary conditioning. For large-scale reservoir simulations involving
tens of millions of cell an exact solution is not feasible, even if
hundreds of processors are used. Nevertheless, the method may be generalised
to approximately solve the system by inverting a preconditioner $B_{i}\thicksim A_{i}$
instead. In particular, ILU(n), nested factorisation, or even multigrid
are all options. Finally, it is also possible to combine boundary
conditioning with a global multigrid pressure solver step \textendash{}
however the use of the far-field means that the solver already includes
some elements of a multigrid solver.

\bibliographystyle{siam}
\bibliography{RoxarReferences}

\pagebreak{}

\appendix

\section*{Appendix}

\section{Error matrix for nested factorisation}

For a hepta-banded matrix of the form

\begin{equation}
A=D+L_{1}+U_{1}+L_{2}+U_{2}+L_{3}+U_{3}
\end{equation}

where $L_{i}$ , $U_{i}$ represent lower and upper bands respectively,
as typically found from three-dimensional finite volume problems,
the nested factorisation preconditioner $B$ is defined by

\begin{equation}
B=\left(P+L_{3}\right)\left(I+P^{-1}U_{3}\right)
\end{equation}

\begin{equation}
P=\left(T+L_{2}\right)\left(I+T^{-1}U_{2}\right)
\end{equation}

\begin{equation}
T=\left(G+L_{1}\right)\left(I+G^{-1}U_{1}\right)
\end{equation}

where $G$ is a diagonal matrix. The error matrix is therefore given
by

\begin{equation}
E=B-A=G-D+L_{1}G^{-1}U_{1}+L_{2}T^{-1}U_{2}+L_{3}P^{-1}U_{3}
\end{equation}

and so we see that by choosing $G$ such that 

\begin{equation}
G=D-L_{1}G^{-1}U_{1}
\end{equation}

we eliminate errors coming from the (strongest) inner bands. In practice,
one may also preserve residual sum by choosing $G$ such that

\begin{equation}
G=D-L_{1}G^{-1}U_{1}-colsum\left(L_{2}T^{-1}U_{2}\right)-colsum\left(L_{3}P^{-1}U_{3}\right)
\end{equation}

which nonetheless preserves the virtue of eliminating errors from
the inner bands.

\section{Motivation for use of lumping}

Suppose we are iteratively solving the $n\times n$ linear system
$r=b-ax=0$. Define the $n\times N$ lumping matrix $E$ consisting
of zeros and unit elements, so that, for any $i$, only one element
$E_{Ii}$ is non-zero. A typical $E$ matrix is shown below

\begin{equation}
E=\left(\begin{array}{ccccccccc}
1 & 1 & 1 & 0 & 0 & \cdots & 0 & 0 & 0\\
0 & 0 & 0 & 1 & 1 & \cdots & 0 & 0 & 0\\
0 & 0 & 0 & 0 & 0 & \cdots & 1 & 1 & 1
\end{array}\right)
\end{equation}

Premultiplying a $n$-dimensional vector by $E$ transforms it to
the $N$-dimensional space (where it this case $N=3$)

\begin{equation}
X=Ex,\qquad X_{I}=\sum_{i\in I}E_{Ii}x_{i}
\end{equation}

A coarse-grid lumping matrix, or Watts correction \cite{Wat71}, can
then be defined as

\begin{equation}
A=EaE^{T},\qquad A_{IJ}=\sum_{i\in I}\sum_{j\in J}E_{Ii}a_{ij}E_{jJ}^{T}
\end{equation}

so that $A_{IJ}$ is the sum of the matrix elements projected out
by row vectors $E_{I}$ and $E_{J}$. Starting with some residual
$r$, one can construct the lumped residual $R=Er$ and then the lumped
search direction is given by

\begin{equation}
\Delta X=A^{-1}R
\end{equation}

Projecting $\Delta X$ back to the $n$-space gives us the search
direction $\Delta x=E^{T}\Delta X$. Then the change in residual is
$\Delta r=-a\Delta x$ and so the new residual is given by

\begin{equation}
r'=r-a\Delta x=r-aE^{T}\Delta X
\end{equation}

Lumping the residual into the $N$-space we find that

\begin{equation}
R'=Er'
\end{equation}

\[
=Er-EaE^{T}\Delta X
\]

\[
=Er-EaE^{T}A^{-1}Er
\]

\[
=Er-AA^{-1}Er
\]

\[
=0
\]

i.e. the residuals that are lumped together by the matrix $E$ sum
to zero.
\end{document}